 \documentclass[11pt,letterpaper]{amsart}
 \usepackage{amsthm, amsxtra, latexsym, mathrsfs, amscd}
 \usepackage{graphicx}
\usepackage{epsfig}
\usepackage{amsfonts}
\usepackage{amscd}
\usepackage{amsfonts}
\usepackage{color}
\usepackage{tikz}
\usepackage{psfrag}
\usepackage{pstricks}
\usepackage{mathdots}
\usepackage{psfrag}
\usepackage{amsthm,amsmath,amssymb,graphics,hyperref}
\usepackage{amsmath,amssymb,hyperref} 
\usepackage{enumerate}
\usepackage[titletoc,title]{appendix}
\input xy
\xyoption{all}
\usepackage{ulem}

\theoremstyle{plain}
\newtheorem{theorem}{Theorem}[section]
\newtheorem{proposition}[theorem]{Proposition}
\newtheorem{lemma}[theorem]{Lemma}
\newtheorem{corollary}[theorem]{Corollary}
\theoremstyle{definition}

\newtheorem{definition}[theorem]{Definition}
\newtheorem{remark}[theorem]{Remark}

\newcommand{\R}{\mathbb{R}}
\newcommand{\bc}{\mathbb{C}}
\newcommand{\ra}{{\rightarrow}}

\let\cal\mathcal
\newcommand{\bR}{{\mathbb R}}
\newcommand{\rsl}{{\mathrm{SL(3, \mathbb R)}}}

\begin{document}
\title [SRB Entropy of  Convex Projective Structure]
      {  Sinai-Ruelle-Bowen measure Entropy of  geodesic flow  on Convex Projective Surfaces  }
        \author{Patrick Foulon and Inkang Kim}

        \date{}
        \maketitle
\begin{abstract}
We study the entropy of  Sinai-Ruelle-Bowen measure of the geodesic flow on convex real projective surfaces, and shows that the Hilbert area tends to infinity  if  the entropy tends to zero. For
 the Blaschke metric, the area tends to infinity if and only if the entropy tends to zero.
\end{abstract}
 \footnotetext[3]{Research by Inkang Kim is partially supported by Grant
NRF-2019R1A2C1083865, KIAS Individual Grant (MG031408) and he also thanks CIRM for the warm support  during his
visit, and IHP for the support through research in Paris program.}

\section{Introduction}
Let $\cal T$ be the Teichm\"uller space of $S$, a closed surface of genus $\geq 2$, $\cal M=\cal T/Mod(S)$ the moduli space, $\cal P$ the space of the marked strictly convex real projective structures, which can be identified as a holomorphic vector bundle over $\cal T$ whose fibres are holomorphic cubic differentials.
Recently real projective structures on surface have drawn much attention as a special but important case of higher Teichm\"uller
theory. People study its geometric aspect, metric aspect and dynamical aspect using different methods and various perspectives.
In this paper, we study its dynamical properties in terms of Sinai-Ruelle-Bowen measure, abbreviated SRB measure, on the  (projectivized) tangent bundle $(TS\setminus \{0\})/\R^*_+$ of a convex real projectrive surface $S$. Until now, quite a few scholars studied its dynamical properties using Bowen-Margulis measure, but we will use
SRB measure to catch  different aspects of the geometry, see for more 
details   \cite{FK2}. 

If the topological entropy tends to zero, then the Hilbert volume of the convex projective manifold tends to infinity \cite{Ad}, but
the converse is not true \cite{FK} in general.  If we replace  Bowen-Margulis measure by SRB measure, the story is quite different. 

Goldman introduced the deformation of projective structures, called bulging deformation along simple closed curves. 
In \cite{FK}, the authors studied its behavior with respect to Bowen-Margulis measure and showed that its topological entropy 
does not go to zero when the parameter of the bulging deformation goes to infinity. On the contrary, it has been expected that
the entropy of such bulging deformation tends to zero with respect to SRB measure \cite{PK}.
In this  article, we show that the volume tends to infinity if  the SRB measure entropy tends to zero with respect to  the Hilbert metric. For the Blaschke metric, we have an if and only if result.

Firstly, we study the entropy with respect to the Hilbert  metric, 
and using refined and sophisticated arguments, we  show
\begin{theorem}\label{main1}Let $S_t$ be a smooth parameter of  deformations of strictly convex real projective structures on a  surface of negative Euler number.
Then the Hilbert area tends to infinity if its measure entropy $h^H_{SRB}(\phi_t)\ra 0$  where $\phi_t$ is the geodesic flow on $HS_t=(TS_t\setminus \{0\})/\R^*_+$ with respect to Hilbert metric.
\end{theorem}

The proof of this theorem is quite involved. Several different areas should contribute to the resolution; ergodic theory, affine and projective geometry, Higgs bundle type theory to analyze the degenerating behavior. 

 First we need to know the relationship between the measure entropy  $h^s_{SRB}$ and the approximately regular exponent $\alpha_s$ along the deformation;
$$\alpha_s=\frac{2}{h^s_{SRB}}.$$ This is carried out in \cite{PK}. Then we have to relate the area to $\alpha_s$. This is done via the $\beta$-convexity as investigated by Benoist in \cite{BIHES}. We analyze the degenerating case via the language of regular separated neck \cite{Loftin} utilizing the holomorphic vector bundle structure over the Teichm\"uller space with fibres being the set of  holomorphic cubic differentials.
We also investigate Gromov's $\delta$-hyperbolicity as in Benoist \cite{BIHES} to relate convexity to other quantities.

The Blaschke metric case is a bit easier to deal with, and we obtain a very satisfying result.
We use  Katok \cite{Katok} and Osserman-Sarnak \cite{OS} inequalities.
\begin{theorem}\label{Blaschke}Let $S_t$ be a continuous parameter of  deformations   of strictly convex real projective structures on a  surface of negative Euler number.
Then its SRB measure entropy $h^B_{SRB}(\phi_t)\ra 0$ as $t\ra\infty$ where $\phi_t$ is a geodesic flow on $HS_t=(TS_t\setminus \{0\})/\R^*_+$ with respect to Blaschke metric if and only if its Hilbert (equiv. Blaschke) area tends to infinity.
\end{theorem}
Note that the above theorem is not true for topological entropy due to the equivalence of Blaschke metric and Hilbert metric.
As we can see in the above theorem, SRB measure better characterizes the relationship between the area and the  measure entropy. In higher dimensional case, the relation still needs to be explored, and the other direction for the Hilbert metric case still remains open.
\section{Preliminaries}
\subsection{Projective structures, Blaschke and Hilbert metrics}\label{Bl}
Let $M$ be an $n$-dimensional closed manifold equipped with a strictly convex projective structure so that
$D:\widetilde M\ra \mathbb{RP}^n$ is a developing map with $D(\widetilde M)=\Omega$. Then it inherits 
a Hilbert metric defined in the following way. Let $x,y\in M$. Then choose $\tilde x,\tilde y$ in the same fundamental domain, and define $d(x,y):=d_\Omega(D(\tilde x), D(\tilde y))$ where the Hilbert metric $d_\Omega$ defined
as below and the corresponding Finsler norm on $T\Omega$ is denoted by $F$. More precisely,
 for $x\neq y\in \Omega$, let $p,q$ be the
intersection points of the line $\overline{xy}$ with $\partial\Omega$ such that $p,x,y,q$ are in this order. The
Hilbert distance is defined by
$$d_\Omega(x,y)=\frac{1}{2}\log \frac{|p-y||q-x|}{|p-x||q-y|}$$ where $| \cdot |$ is a Euclidean
norm in an affine chart containing $\Omega$. This metric coincides with the hyperbolic metric if
$\partial\Omega$ is a conic. Note that the Hilbert metric on $M$ also depends on the boundary of $\Omega$. The Hilbert metric is Finsler rather than
Riemannian. The Finsler norm $F=||\cdot||$ is given, for $x\in\Omega$ and a
vector $v$ at $x$, by
$$||v||_x=\frac{1}{2}\big(\frac{1}{|x-p^-|}+\frac{1}{|x-p^+|}\big)|v|$$ where $p^\pm$
are the intersection points of the line with $\partial\Omega$, defined by $x$ and $v$ with the
obvious orientation, and where $| \cdot |$ is again 
a Euclidean norm.

Another way to describe projective structures is by using hypersurface theory.
For a hypersurface $\mathcal H\subset \R^{n+1}$ with a trivial bundle $E=\mathcal H\times \R^{n+1}$, choose a transverse vector field $\xi$ over $\mathcal H$ so that $E=T\mathcal H\oplus \cal L$ where $\cal L$ is a trivial line bundle over $\mathcal H$ spanned by $\xi$. For the standard affine flat
connection $\nabla$  on $\R^{n+1}$, its restriction on $E$ gives
\begin{align}\label{affine metric}
 &\nabla_X Y=D_XY + h(X,Y)\xi\in T\mathcal H\oplus \cal L,\\
 &\nabla_X \xi=-S(X)+ \tau(X)\xi \in T\mathcal H\oplus \cal L,
 \end{align}
where $D$ is a torsion-free connection on $T\mathcal H$, $h$ is a symmetric 2-form on $T\mathcal H$, $S$ is a shape operator and $\tau$ is a 1-form. If $\mathcal H$ is locally strictly convex, i.e., it can be written locally as the graph of a function with positive definite hessian, then there exists a unique $\xi$ such that
\begin{itemize}
\item $\tau=0$, $h$ is positive definite and
 \item $|det(Y_1,\cdots,Y_n,\xi)|=1$ for any $h$-orthonormal frame $Y_i$ of $T\mathcal H$.
 \end{itemize}
 This vector field $\xi$ is called an affine normal, $D$ Blaschke connection, $h$ affine or Blaschke metric. These are invariant under $SL^\pm(n+1,\mathbb R)$ of real matrices with determinant $\pm 1$.
  If the shape operator $S=-Id$, the hypersurface $\mathcal H$ is an affine sphere with affine curvature $-1$.
  For such an affine sphere, $O=x-\xi(x)$ is a point, called the center of the affine sphere.
  Cheng and Yau \cite{CY} solved the following conjecture of E. Calabi \cite{Cal}.
\begin{theorem}If $\cal C\subset\R^{n+1}$ is an open convex cone containing no lines,
there exists a unique embedded affine sphere $\mathcal H$ whose center is origin and affine curvature $-1$,
which is asymptotic to the boundary of $\mathcal C$.
\end{theorem}
If $\Omega\subset\R^n$ is a bounded open convex set so that
$$\cal C=\{t(1,x)| x\in\Omega, t>0\}$$ then
$$\cal H=\{\frac{-1}{u(x)}(1,x)| x\in\Omega\}$$ where $u$ is the unique convex solution
of the real Monge-Amp\`ere equation
$$det D^2 u=(-1/u)^{n+2}$$ under the condition that
$u|_{\partial \Omega}=0.$
Using a Benz\'ecri's theorem \cite{Be},
Benoist-Hulin showed \cite{BH1}[Prop 3.4] that there exists a positive constant $C_n$ depending only on the dimension such that for any tangent vector $X$,
\begin{equation}\label{comp}
\frac{1}{C_n} ||X||_H\leq ||X||_B \leq C_n ||X||_H,
\end{equation}
 where $H$ and $B$ stand for Hilbert and Blaschke respectively.

For strictly convex domain $\Omega$, we consider a Pick form $c$ which is defined as
$$h((D_X- \nabla^h_X)Y, Z)=c(X,Y,Z),$$ where $\nabla^h$ is the Riemannian connection defined by the Blaschke metric $h$.
This Pick form is known to be a real part of a cubic holomorphic differential  with respect to the underlying conformal hyperbolic metric on the projective surface. Using these data,  Labourie \cite{La} and Loftin \cite{Loftin1} independently showed that the space of strictly convex real projective structures on a closed surface of genus at least 2 is a holomoprhic vector bundle over Teichm\"uller space whose fibers are cubic holomoprhic differentials.

If $h$ is the Blaschke metric determined by the cubic differential $b=f(z)dz^3$ with respect to the conformal structure underlying $h$, 
Wang's equation states that the curvature of $h$ is 
\begin{eqnarray}\label{curv}
-1\leq \kappa_{h}=-1+2||b||^2_{h} < 0,\end{eqnarray} where $||b||^2_{h}=e^{-3(\psi+\phi)}|f(z)|^2$ is the pointwise norm of $b$ with respect to $h=e^{\psi}e^{\phi}|dz|^2$ with $e^\psi|dz|^2$ being the underlying hyperbolic metric.
See \cite{BH} for strict negativeness of the curvature.

\subsection{Lyapunov exponents and approximately regular value}
Let $\phi=\phi^t$ be a $C^1$ flow on a Riemannian manifold $W$. A point $w\in W$ is said to be regular if there exists
a $\phi^t$-invariant decomposition along $\phi^t w$
$$TW=E_1\oplus \cdots \oplus E_p$$  and real numbers
$$\chi_1(w)<\cdots< \chi_p(w),$$ such that, for any vector $Z_i\in E_i\setminus \{0\}$,
$$\lim_{t\ra\pm \infty}\frac{1}{t}\log ||d\phi^t(Z_i)||=\chi_i(w),$$ and
\begin{eqnarray}\label{det}
\lim_{t\ra\pm \infty}\frac{1}{t}\log |\text{det} d\phi^t|=\sum_{i=1}^p\text{dim}E_i \cdot \chi_i(w).
\end{eqnarray}
The numbers $\chi_i(w)$ associated with a regular point $w$ are called the Lyapunov exponents of the flow at $w$.
Due to Oseledec's  multiplicative ergodic theorem \cite{Os}, the set of regular points has  full measure with respect to any flow invariant probability measure.

For 2-dimensional convex domain $\Omega$ with $C^1$-boundary, the neighborhood of each point of the boundary $\partial \Omega$ can be written as the graph of a convex function $f$ near the origin. Following Crampon \cite{Cr}, such a function $f$ is said to be {\it approximately $\alpha$-regular} at the origin for an $\alpha\in [1,\infty]$, if
$$\lim_{t\ra 0} \frac{\log \frac{f(t)+f(-t)}{2}}{\log |t|}=\alpha.$$ This quantity is invariant under affine and projective
transformations. When $f$ is quasi-symmetric, approximately $\alpha$-regular for $\alpha$ finite means that the function behaves like $|t|^\alpha$ near the origin. 

In general,  even though $\alpha$ exists for the set of full Lebesgue measure on $\partial\Omega$ for the convex domain $\Omega$ whose quotient has a finite volume,  $\alpha(x)$ might not exist for some points on $\partial \Omega$, hence we define
$$\alpha_{sup}(x)=\limsup_{t\ra 0} \frac{\log \frac{f(t)+f(-t)}{2}}{\log |t|},$$ for a convex function
$f$ whose graph around the origin is the neighborhood of $x$ on $\partial \Omega$.
For a hyperbolic isometry $\gamma$ with eigenvalues $\lambda_1>\lambda_2>\lambda_3$, 
$$\alpha(\gamma^+)=\frac{\log \lambda_1/\lambda_3}{\log\lambda_1/\lambda_2}.$$

\subsection{SRB measure} 
It is known that the Bowen-Margulis measure maximizes the entropy among the geodesic flow invariant probability measures \cite{crampon}. This entropy is equal to the topological entropy $h_{top}(\phi)$  of the geodesic flow $\phi$ \cite{crampon} and it is also equal to the exponential growth of the lengths of  closed geodesics:
$$h_{top}(\phi)=\lim_{R\ra\infty} \frac{\log \#\{[\gamma]|\ell(\gamma)\leq R\}}{R}.$$ For strictly convex real projective case, this topological entropy is also equal to the volume entropy of the Hilbert metric \cite{crampon}.  A general result of Crampon \cite{crampon} implies that the entropy of Bowen-Margulis measure of strictly convex real projective structures  on a closed surface lies in $(0,1]$ and it is 1 if and only if the Hilbert metric is hyperbolic.

There exists another important invariant probability measure whose existence is proven at least for Anosov flows. The following are  some characterizations of SRB measure.
\begin{enumerate}
\item  The measure attains the equality in Ruelle inequality $$h_\mu(\phi)\leq \int \chi^+ d\mu,$$ where $\chi^+=\sum \text{dim} E_i\cdot \chi_i^+$ denotes the sum of positive Lyapunov exponents.
\item  There exists a set $V$ of full Lebesgue measure
such that for each continuous function $f:M\ra \R$ and for every $x\in V$,
$$\lim_{T\ra\infty} \frac{1}{T}\int_0^T f(\phi^s(x))ds=\int f d\mu_{SRB}.$$
\item   The measure has Lebesgue absolutely continuous conditional measures on unstable manifolds. This characterization is due to 
Ledrappier-Young \cite{LY}.
\end{enumerate}

For hyperbolic surface, Bowen-Margulis measure equals SRB measure, but in general they are quite different. If there exists an invariant normalized volume form $\omega$, then it coincides with SRB measure.  Even  when volume is not preserved, SRB measures are somehow the invariant measures most compatible with volume. They provide a mechanism for explaining how local instability on attractors can produce coherent statistics for orbits starting from large sets in the basin. See the survey article \cite{Young}.

In \cite{PK}, it is shown that $h^s_{SRB}=1+\eta_s$ along the smooth deformation where $h^s_{SRB}$ is the SRB measure entropy and $\eta_s\leq 0$ is the SRB almost sure value of the sum of the parallel Lyapunov exponents. Furthermore, there exists a subset of full Lebesgue measure $\cal E_s\subset \partial \Omega_s$ such that for any $\psi\in \cal E$, the approximately regular exponent
$$ \alpha_s=\alpha(\psi)=\frac{2}{1+\eta_s}\geq 2.$$
In addition, $\eta_s$ and $\alpha_s$ are continuous functions in $s$.
Indeed, the denominator $1+\eta_s$ is equal to $\chi^+_s$, the sum of positive Lyapunov exponents, which is again a continuous function in $s$. 
Note that $\alpha_s\ra \infty$ if $h^s_{SRB}\ra 0$.

We will introduce $\beta$-convexity, which is bigger than $\alpha$, in the next section. These quantities $\beta$ and $\alpha$ will be essential to the proof of our main theorem.

\section{$\beta$-convexity, $\delta$-hyperbolicity and $H$-quasisymmetricity}
The content in this section is based on \cite{BIHES}. Readers can skip this section for the first reading since we just use theorems listed out in this section. In our case, we  work on $\Omega\subset \mathbb{RP}^m$, and $\partial\Omega$ which can be locally represented by the graph of a strictly convex $C^1$ function $f:U \ra \mathbb R$ for some  convex open set $U\subset \mathbb{R}^{m-1}$.

This section is technical, but important in affine and projective geometry. It deals with
 Gromov's hyperbolicity of the domain in relation to $\beta$-convexity. This relation will be crucial for the proof of our main theorem \ref{main1}. We outline the argument, rather than giving a whole spectrucm of the theory.
 First we begin with the definition of $\beta$-convexity.
\begin{definition}[Def 1.2 in \cite{BIHES}] For a $C^1$-convex function $F:U\ra \mathbb R$, set $D^F_x(h)=F(x+h)-F(x)- F'(x)h$.  It is said to be  quasisymmetrically convex if there exists $H\geq 1$ such that
$$D^F_x(h)\leq H D^F_x(-h),$$ whenever $x-h$ and $x+h$ are in $U\subset \mathbb R^{m-1}$.

For $\beta\in (2,\infty)$, $F$ is said to be {\bf $\beta$-convex} if $$\inf_{\{(x,h):h\neq 0\}} |h|^{-\beta}D^F_x(h)>0.$$  
\end{definition}
For a given strictly convex   domain $\Omega\subset \mathbb{PR}^m$, one can cover $\partial \Omega$ by small open sets $\cal O_i$ and each $\cal O_i$ is realized by the graph of a strictly convex $C^1$ function $F_i:U_i \ra \mathbb R$ where $U_i$ is a convex open set in $\mathbb R^{m-1}$.
Then set $\beta(\Omega)$ to be $$\max_i\min(\beta: F_i\ \text{is}\ \beta\ \text{convex on all compact sets of}\ \cal O_i)$$
It is proved in Corollary 1.5 (b) of \cite{BIHES} that if $\Omega$ is (Gromov) hyperbolic, then $\partial\Omega$ is $\beta$-convex for some $\beta\in (2,\infty)$.

We list out some theorems in \cite{BIHES} which are crucial to our paper.
First we fix some notations.

Let $X_m$ be  the set of properly convex open domains in $\mathbb{RP}^m$ with Hausdorff Topology. It is known \cite{Be} that $G_m=PGL(m+1,\mathbb R)$ acts properly cocompactly on $X_{m,0}=\{(\Omega,x):\Omega\in X_m, x\in\Omega\}$ (Prop 2.3 \cite{BIHES}). Denote $X_m^\delta\subset X_m$ the set of strictly convex $\delta$-hyperbolic open domains.

Benoist \cite{BIHES} showed the following properties.
\begin{proposition}[\cite{BIHES}, Prop. 2.10]
\begin{enumerate}[(a)]
\item For $\delta>0$, $X_m^\delta$ is a closed $G_m$-invariant set whose every element is strictly convex.
\item If $F$ is a closed $G_m$-invariant subset of $X_m$, whose every element is strictly convex, then $F\subset X^\delta_m$ for some $\delta>0$.
\end{enumerate}
\end{proposition}
\begin{corollary}[\cite{BIHES}, Cor. 1.5]\label{deltahyperbolic} For a properly convex open domain $\Omega\subset \mathbb{RP}^m$,
\begin{enumerate}[(a)]
\item If $\Omega$ is $\delta$-hyperbolic, then $\partial \Omega$ is $C^\alpha$ for $\alpha\in (1,2)$.
\item If $\Omega$ is $\delta$-hyperbolic, then $\partial \Omega$ is $\beta$-convex for $\beta\in (2,\infty)$.
\item When $\partial\Omega$ is real analytic, $\Omega$ is $\delta$-hyperbolic.
\end{enumerate}
\end{corollary}
\begin{corollary}[\cite{BIHES}, Cor. 2.9] Suppose $\Omega$ is a properly convex open domain in $\mathbb{ RP}^m$. Suppose $\Omega$ is not strictly convex, or $\partial \Omega$ is not $C^1$. Then 
the closure of $G_m\Omega$ contains a properly convex open domain
 which admits a triangular section.
\end{corollary}

Guichard \cite{Gu} showed that $\beta(\Omega)=\sup_\gamma \beta(\gamma)$ for a strictly convex divisible domain $\Omega$, and record
 the following lemma.
\begin{lemma}\label{continuity}$\beta:\{\Omega: \Omega\ \text{  strictly convex divisible}\}\ra (2,\infty)$ is a lower-semicontinuous function with respect to the Gromov-Hausdorff topology.
\end{lemma}
\begin{proof} 
The $\beta$ of $\gamma$ is given by
$$ \beta(\gamma)=\frac{\log \lambda_1-\log\lambda_{3}}{\log \lambda_1-\log\lambda_{2}}.$$
Since $\beta(\Omega)=\sup_\gamma \beta(\gamma)$ for a strictly convex divisible domain $\Omega$,   $\beta$ function
$$\beta:\{\Omega: \text{strictly convex divisible}\}\ra \mathbb R^+$$ is lower-semicontinuous since $\beta_\gamma(\rho_s)=\beta(\rho_s(\gamma))$ is a continuous function on the set of holonomy representations of strictly convex real projective structures. Since the developing map is a local homeomorphism
from the set of holonomy representations, up to conjugacy and up to mapping class group action, to the developed images of the projective structures up to the action of $G_2$ (see Prop. 2.6 in \cite{PK}), we proved the claim.
\end{proof}

In the following lemma, we deal with a finite volume strictly convex projective surface
to relate two quantities, $\alpha$ and $\beta$.
\begin{lemma}\label{lem:3.2} If $\Omega$ admits a quotient which is a finite volume strictly convex real projective surface, then 
$1<\sup_{x\in\partial\Omega}\alpha_{sup}(x)\leq \beta(\Omega)<\infty$ where $\beta(\Omega)$ is the beta convexity for $\Omega$.
\end{lemma} 
\begin{proof}Benoist-Hulin  showed that the curvature of Blaschke metric is negatively curved (Proposition 3.3 in \cite{BH})
and approaches a negative constant deep into a cusp (Proposition 3.1 in \cite{BH1}), hence the curvature is pinched negative on the surface. Since the Hilbert metric and the Blaschke metric are comparable \cite[Corollary 4.7]{BH}, $\Omega$ equipped with the Hilbert metric is Gromov hyperbolic. Furthermore Benoist  showed that $\Omega$ is quasisymmetrically convex, hence $\partial \Omega$ is $\beta$-convex for some $\beta\in (2,\infty)$ (Corollary 1.5 (b) in \cite{BIHES}).

Hence for each $x\in \partial \Omega$, $\partial\Omega$ near $x$ can be represented by the graph of a quasisymmetrically  $\beta$-convex $C^{1}$ function $f$ with $f(x)=f'(x)=0$. 

Let  $\alpha_{sup}(x)=\alpha$  for  $x\in \partial\Omega$.  Then it is easy to see that for any $\epsilon>0$, there exists a sequence $t_n\ra 0$ (see also Lemma 3.4.3, \cite{Cr}), such that $$ \frac{f(x+t_n)+f(x-t_n)}{2} \leq |t_n|^{\alpha-\epsilon}.$$  Hence $$ f(x+t_n)+f(x-t_n)= D_x(t_n)+ D_x(-t_n)\leq 2|t_n|^{\alpha-\epsilon}.$$ By the $\beta$-convexity of $\partial\Omega$, $D_x(-t)> c |t|^\beta$ for small $t$ and for some fixed constant $c>0$. Consequently, for any small $\epsilon>0$, there exists a sequence $t_n\ra 0$ such that
$$2c|t_n|^{\beta}\leq D_x(t_n)+ D_x(-t_n) \leq 2|t_n|^{\alpha-\epsilon} .$$ 
This is possible only when $\alpha\leq \beta$. Hence the claim follows.
\end{proof}
From this lemma, we know that $\alpha(\Omega)\leq \beta(\Omega)$ where $\alpha(\Omega)$ is the almost sure value defined on the full Lebesgue measure
$\cal E\subset \partial\Omega$ satisfying   $\alpha(\Omega)=\frac{2}{h_{SRB}}$.

\section{Holomorphic cubic differentials and degeneration of convex projective structures on surface}\label{neck}
\subsection{Cubic holomorphic differentials}
Let $\Omega\subset \bc$ be a convex bounded open domain such that its
corresponding hyperbolic affine sphere $\cal H$ in $\R^3$ can be described as
a conformal map  $\xi:\Omega\ra \cal H$ with respect
to the standard conformal structure on $\Omega$ and the conformal
structure induced from the affine metric on $\cal H$. By choosing a
local complex coordinates $z$ on $\Omega$, the structure equation (\ref{affine metric})
of Blaschke connection and Blaschke (affine) metric give rise to
\begin{eqnarray}\label{cubic}
\triangle \phi + 4 e^{-2\phi}||U||^2-2e^\phi+2=0,
\end{eqnarray} where the affine metric
is $e^\phi g$, $g=e^\psi|dz|^2$ is the hyperbolic metric on $\Omega$ and
$U$ is a holomorphic cubic differential. Here
$\triangle=4e^{-\psi}\partial_z\partial_{\bar z}$  is the Laplacian and
$||\cdot||^2=|\cdot|^2e^{-3\psi}$ is the induced
norm from $g$. Note that the cubic differential $U=0$ implies
$\phi=0$, and in this case, the affine structure is just a
Riemannian hyperbolic structure. Again by Benz\'ecri's compactness
theorem, with the notations from Section \ref{Bl}, if we write $U(z)=f(z)dz^3$,
\begin{proposition}[Benoist-Hulin, \cite{BH1}]\label{pick}There exists a universal constant $C$ independent of the projective structure such that the ratio $|U|^{2/3}/\mu_h$  between the affine measure
$\mu_h=(-u)^{-3}|dz|^2$ and the Pick measure $|U|^{\frac{2}{3}}=|f(z)|^{2/3}|dz|^2$ is uniformly bounded by $C$.
\end{proposition}
\begin{lemma}\label{cubic} If $(U_n, g)$ is a sequence of real projective structures whose Blaschke metrics are given by $e^{\phi_n}g$.  If  $||U_n||^2\ra\infty$, then the corresponding Hilbert metric area tends to 
infinity.
\end{lemma}
\begin{proof} If $||U_n||^2\ra\infty$, by \cite{Xin}, the volume entropy of the Hilbert metric tends to zero. Theny by \cite{Ad}[Theorem 1.4], the corresponding Hilbert area tends to infinity.
\end{proof}

\subsection{Degeneration of projective structures}\label{degen}
If $\Omega$ is not a triangle, the elements of $\rsl$ acting on  a properly convex domain $\Omega$ are classified as follows (see \cite{Marquis}):
\begin{enumerate}
\item {\bf Hyperbolic :} the matrix is conjuagate to $$\begin{bmatrix} 
           \lambda^+ & 0 & 0 \\
            0              & \lambda^0 & 0\\
            0           & 0 & \lambda^-\end{bmatrix}, \begin{array}{l}  \text{where } \lambda^+>\lambda^0>\lambda^->0 \\ \text{and } \lambda^+\lambda^0\lambda^-=1.\end{array}$$
\item {\bf Quasi-hyperbolic :} the matrix is conjugate to $$\begin{bmatrix} 
           \alpha & 1 & 0 \\
            0              & \alpha & 0\\
            0           & 0 & \beta \end{bmatrix},  \begin{array}{l}  \text{where } \alpha, \beta>0, \alpha^2 \beta=1 \\ \text{and }\alpha, \beta \neq 1.\end{array}$$
\item {\bf Parabolic :}  the matrix is conjugate to  $$\begin{bmatrix} 
           1 & 1 & 0 \\
            0              & 1 & 1\\
            0           & 0 & 1\end{bmatrix}.$$
\item {\bf Elliptic :}  the matrix is conjugate to $$\begin{bmatrix} 
           1 & 0 & 0 \\
            0              & \cos\theta & -\sin\theta \\
            0           & \sin\theta & \cos\theta\end{bmatrix},\ \text{where } 0<\theta<2\pi.$$
\end{enumerate}

Choose coordinates in $\mathbb{RP}^2$ so that the
hyperbolic action is given by the diagonal action $(\lambda_1,\lambda_2,\lambda_3)$. The three fixed points
of this action are the attracting fixed point [1, 0, 0], the repelling fixed
point [0, 0, 1], and the saddle fixed point [0, 1, 0]. Define the {\bf principal
triangle} T as the projection onto $\mathbb{RP}^2$ of the first octant in $\mathbb{R}^3$. The quotient of this principal triangle by the hyperbolic isometry is called the {\bf principal half-annulus}, which has an infinite length (and an infinite volume) with respect to the Hilbert metric. The
{\bf principal geodesic} $\tilde l$ associated to this  matrix is the straight
line in the boundary of T from the repelling to the attracting fixed
point.

The Hilbert length of the action of the hyperbolic isometry $\gamma$ is given by
$$\ell(\gamma)=\frac{1}{2}\log \frac{\lambda_1}{\lambda_3}.$$

\smallskip
If one pinches down a simple closed geodesic in a hyperbolic surface, the hyperbolic surface will converge to a hyperbolic surface with a neck pinched along that closed geodesic in Mumford compactification of the Teichm\"uller space.
Such a phenomenon occurs even in the real projective stuctures even though the Hilbert length along the neck is not necessarily zero.

An end of a convex $\mathbb{RP}^2$ surface with hyperbolic holonomy has {\bf bulge
$\infty$} if the image of the developing map contains a principal triangle T.
Similarly, an end of a convex $\mathbb{RP}^2$ surface with hyperbolic holonomy
has {\bf bulge $-\infty$} if the image 
 of the developing map has the principal
geodesic $\tilde l$ in its boundary.  See the statement before Section 2.4 in \cite{Loftin}.

A pair of ends of a convex projective surface forms a {\bf regular separated neck} in three categories: See the statement before Theorem 2.6.1  and Theorem 4.1.1 in \cite{Loftin}.
\begin{enumerate}
\item The holonomy around each end is parabolic
\item The holonomy around each end is quasi-hyperbolic,
\item The holonomy around each end is hyperbolic; the oriented holonomies
around each end are, up to conjugation, inverses of each other;
and the $\mathbb{RP}^2$ structure about one of the two ends has bulge $\infty$,
while the other end has bulge $-\infty$.
\end{enumerate}
Note that in the first case, the Hilbert length of the holonomy is zero, and  for the last two cases, the Hilbert length of the holonomy is positive.

A simple end of a convex projective surface is  called regular if it forms half of a regular separated neck.
It is shown in \cite{Loftin} that all the neck-separating degenerations of convex projective structure on a surface fall into
the above three categories.
In the latter two cases,   in the universal cover $\Omega$, {\bf either the axis of (quasi)-hyperbolic isometry is a part of the boundary $\partial\Omega$, or the principal triangle is attached along the principal geodesic $\tilde l$ depending on whether the bulge is $-\infty$ or $\infty$.}

%
Note that if a sequence of projective structures degenerate to the one with a (quasi)-hyperbolic regular separated neck $l$, the Hilbert length of any curve transversal to $l$ must go to infinity since the neck will be separated in the limit.
Note also that the corresponding domains $\Omega_s$ converge to the limit domain $\Omega_\infty$ in the Gromov-Haudorff topology under degeneration and we have the following corollary.
\begin{corollary}\label{infinitetube}Let $S_t$ be a smooth deformation of strictly convex projective structures on $S$, which develops a regular separated neck along a simple closed curve $c$ with (quasi)-hyperbolic holonomy at the limit structure. Then the Hilbert area of $S_t$ tends to infinity.
\end{corollary}
\begin{proof}By the classification of degeneration to a regular separated necks, it corresponds to the case (2) and (3) whose holonomy around the neck is (quasi)-hyperbolic.
In this case, in the limit, a principal half-annulus is attached along the boundary on one end whose bulging parameter is $\infty$. The area of this half-annulus is infinity. Hence the
Hilbert area of $S_t$ tends to infinity.
\end{proof}

\section{Hilbert metric case}

In this section, we show our first main theorem;  $h_{SRB}^s\ra 0 \Rightarrow  \text{area}\ra \infty$.\\

Since $h^s_{SRB}\ra 0$,  $\beta_s\ra\infty$  by Lemma \ref{lem:3.2}, using the fact that $\alpha_s=\frac{2}{h^s_{SRB}}\leq \beta_s$.
Recall that $\cal T$ is the Teichm\"uller space of $S$, $\cal M=\cal T/Mod(S)$ the moduli space, $\cal P$ the space of the marked strictly convex projective structures, which can be identified with a holomorphic vector bundle over $\cal T$ whose fibres are holomorphic cubic differentials.

Let $S_s$ be such a smooth family of deformations of projective structures with $h^s_{SRB}\ra 0$. Let $S_s$ correspond to $(g_s, U_s)$ where
$U_s$ is a holomorphic cubic differentials on the hyperbolic surface with metric $g_s$.

{\bf We show that $(g_s,U_s)$ diverges in a way that either the vertical direction $U_s$ diverges, or the sequence develops a  regular separated neck in the limit whose holonomy is (quai)-hyperbolic.}
\vskip .1 in
CASE i)  ${(g_s, U_s)}$ is contained in a compact set of $\cal P$.

In this case there is a limit  $S_s\ra S_\infty \in \cal P$ after passing to a subsequence. Then 
by the continuity of the entropy map,  $h^s_{SRB}\ra h^\infty_{SRB}>0$, which is a contradiction to the assumption that $h^s_{SRB}\ra 0$.
\vskip 0.1 in
CASE ii)    ${(g_s, U_s)}$ is not contained in a compact set of $\cal P$.\\
There are two cases to consider in this case.

$\bullet$ Suppose the image of $g_s$ lie in a compact set of the moduli space $\cal M$ of $S$. By applying elements in the mapping class group to $S_s$, we may assume that
the set $\{g_s\}$ lie in a compact set of Teichm\"uller space. 
Hence  $||U_s||^2\ra\infty$.  Then the Hilbert area of $S_s$ must diverge to $\infty$ by Lemma \ref{cubic}. In this case, even the topological entropy tends to zero by \cite{Xin}.

$\bullet$ Suppose that $g_s$ does not lie in a compact set of moduli space $\cal M$. By Mumford compactness theorem, the systole length of $g_s$ must go to zero.  Since the action of the mapping class group on the curve complex is cocompact (even on the pants complex),
by applying mapping class group again and passing to a subsequence, we can assume that the hyperbolic length $\ell_{g_s}(c)\ra 0$ for any $c$ in a finite system $C$ of disjoint, non-homotopic
simple closed curves  on $S$. Then the degenerating projective structure is regular as described in Section \ref{degen}, see Prop. 2.9.2 in \cite{Loftin}. Let $\Omega_s\ra \Omega$ in the Gromov-Hausdorff topology. There are two subcases to consider in this setting.

SUBCASE I) Suppose  all of the Hilbert lengths of $c$ in $C$ tend to zero.

In this case, all the necks become  parabolic, and the corresponding $\Omega$ is strictly convex  admitting a finite volume quotient. Hence it is $\delta$-hyperbolic by the argument in the proof of Lemma \ref{lem:3.2}.  Thus $\Omega$ is $\beta$ convex for a finite $\beta>0$ by Corollary 1.5 (b) in \cite{BIHES}
 (Corollary     \ref{deltahyperbolic}), which is a contradiction to $\beta_s\ra \infty$.

SUBCASE II)  Some of the Hilbert length of $c$ in $C$ tend to positive numbers.

These correspond to   last  two cases of regular degenerating cases, and   the area of at least one end  tends to infinity  by Corollary \ref{infinitetube}.

In conclusion, the area will tend to infinity if the SRB measure entropy tends to zero. This finishes the proof.\\

\section{Blaschke metric case}
Let $(M, \sigma)$ be a Riemann surface and $h^\lambda_\sigma$ the measure entropy with respect to the Liouville measure
on the unit tangent bundle. The total area of the surface is denoted by $V_\sigma$ and $d\mu_\sigma=\frac{dvol_\sigma}{V_\sigma}$ the normalized volume form.
Let $\rho$ be a function coming from the conformal equivalence theorem, i.e., $\rho\sigma$ has constant negative curvature with the same volume as $\sigma$, such that $\int \rho d\mu_\sigma=1$.
We set $\int \rho^{1/2} d\mu_\sigma=\rho_\sigma\leq 1$. Then it is proved in \cite{Katok}
\begin{theorem}If $\sigma$ is a metric without focal points, then
$$h_\sigma^\lambda\leq  \rho_\sigma(-2\pi \chi(M)/V_\sigma)^{1/2}.$$
\end{theorem}

If the Hilbert area of projective surface $S_s$ with respect to Hilbert metric  tends to infinity,
by Benoist-Hulin, the area of $S_s$ with respect to Blaschke metric tends to infinity, hence by the previous theorem,
the metric entropy tends to zero.

{\bf Remark} The geodesic flow of the Blaschke  metric preserves a volume form. Hence after normalization, it becomes a geodesic flow invariant probability measure. Then the conditional measure on the unstable foliation is clearly Lebesgue, hence  by the characterization of SRB measure, this normalized Liouville measure is the SRB measure since the Blaschke metric is negatively curved.\\

In the proof above, we only used the fact that the area of the surface with respect to strictly negatively curved Riemannian metrics
tends to infinity. Hence we have
\begin{theorem}The  entropy of SRB measures of  Riemann surface  with curvature $-1\leq \kappa<0$, tends to zero if and only if the area of surface  tends to infinity.  Especially, this holds for the strictly convex real projective surfaces with respect to Blaschke metrics whose Hilbert area tends to infinity.
\end{theorem}
\begin{proof}
Suppose that the SRB entropy tends to zero. 
By Osserman-Sarnak \cite{OS}, if $h_\sigma^\lambda\ra 0$, then the following integral tends to zero.
$$h_\sigma^\lambda\geq \int \sqrt{-\kappa} d\mu_\sigma\geq \int {-\kappa} d\mu_\sigma.$$
Then by Gauss-Bonnet theorem
$$\frac{\int -\kappa dvol_\sigma}{V_\sigma}=\frac{-2\pi\chi(S)}{V_\sigma},$$ hence
$V_\sigma\ra\infty$.
\end{proof}
In \cite{Katok}, it is also shown that the asymptotic Cheeger isoperimetric constant $C_\sigma$, which is defined as the lower limit of the ratios of the length of a rectifiable closed Jordan curve on $\tilde M$ to the area bounded by the curve as the area goes to infinity, satisfies
$$C_\sigma\leq \rho_\sigma( -2\pi \chi(S)/V_\sigma)^{1/2}.$$
Hence as a corollary we get
\begin{corollary}
The asymptotic Cheeger isoperimetric constant with respect to both Hilbert and Blaschke metrics goes to zero as the Hilbert area tends to infinity.
\end{corollary}
\begin{proof}
Since  Blaschke metric and Hilbert metric are comparable, the above inequality holds for both metrics.
\end{proof}

We also have
\begin{corollary}
  SRB measure entropy of Blaschke metric tends to zero if SRB measure entropy tends to zero for convex projective surface with Hilbert metric.
\end{corollary}
\begin{proof} If $h_{SRB}\ra 0$ for Hilbert metric, by Theorem \ref{main1}, the Hilbert area tends to infinity, then Blaschke area also tends to infinity. Then again by Theorem \ref{Blaschke}, the SRB measure entropy with respect to the Blaschke metric tends to zero. 
\end{proof}

\vskip .1 in
\noindent     Patrick Foulon\\ Aix-Marseille Universit\'e, CNRS,  Marseille, France.
\\UMR 822, 163 avenue de Luminy,
13288 Marseille cedex 9, France\\
\texttt{foulon\char`\@cirm-math.fr}\\
\vskip .005 in
\noindent     Inkang Kim\\
     School of Mathematics\\
     KIAS, Heogiro 85, Dongdaemen-gu,
     Seoul, 02455, Korea\\
     \texttt{inkang\char`\@ kias.re.kr}
\end{document}